\documentstyle[amssymb,12pt]{article}
\title{Mutually algebraic structures and expansions by predicates}

\author{Michael C.\
Laskowski\thanks{Partially supported
by NSF grants DMS-0600217 and DMS-0901336.}\\Department of 
Mathematics\\University of Maryland}

\newbox\smilebox
\newbox\anchorbox
\newbox\noanchorbox
\newbox\tempbox

\setbox\smilebox=\hbox{$\smile$}

\def\anchor{\hbox{\vtop{
           \hbox to \wd\smilebox{\hfil\vrule width.4pt height7pt depth1pt\hfil}
           \vskip  -11.5truept
           \hbox to \wd\smilebox{\hfil$\smile$\hfil}}}}
\setbox\anchorbox=\anchor
\def\noanchor{\hbox{\vtop{
           \hbox to \wd\anchorbox{\hfil\anchor\hfil}
           \vskip -14truept
           \hbox to \wd\anchorbox{\hfil/\hfil}}}}
\setbox\noanchorbox=\noanchor

\def\fg#1#2#3{\setbox\tempbox=\hbox{$\scriptstyle{#2}$}
\ifnum\wd\anchorbox>\wd\tempbox\dimen255=\wd\anchorbox
\else\dimen255=\wd\tempbox\fi
{#1\,\vtop{\hbox to \dimen255{\hfil\anchor\hfil}
           \vskip -6truept
           \hbox to \dimen255{\hfil$\scriptstyle{#2}$\hfil}}
           \,#3}}

\def\nfg#1#2#3{\setbox\tempbox=\hbox{$\scriptstyle{#2}$}
\ifnum\wd\noanchorbox>\wd\tempbox\dimen255=\wd\noanchorbox
\else\dimen255=\wd\tempbox\fi
{#1\,\vtop{\hbox to \dimen255{\hfil\noanchor\hfil}
           \vskip -6truept
           \hbox to \dimen255{\hfil$\scriptstyle{#2}$\hfil}}
           \,#3}}

\setbox1=\hbox{$\bot$}

\def\north#1#2{#1\,
\hbox{$\bot$\llap {\hbox to\wd1 {\hfil $/$\hfil}}}
\,#2}

\def\nao#1#2#3{#1\  \hbox{\vtop{ 
\baselineskip=4pt
\hbox{$\bot$\llap {\hbox to\wd1 {\hfil $/$\hfil}}
\hskip .05em \llap{\hbox{$^{\scriptscriptstyle{a}}$}}}\hbox{$\scriptstyle
{#2}$}}}\, #3}

\def\bp{\par{\bf Proof.}$\ \ $}

\def\includeE#1{{\lhook\kern-3.5pt\joinrel\smash{
    \mathop{\longrightarrow}\limits^{#1}}}}

\def\efor/{Example~\ref{E4}}

\def\BL/{Baldwin--Lachlan}
\def\Bu/{Buechler}
\def\Hr/{Hrushovski}
\def\lm/{locally modular}
\def\wm/{weakly minimal}
\def\nm/{non--modular}
\def\ss/{superstable}
\def\ud/{unidimensional}
\def\sm/{strongly minimal}

\def\abar{\overline{a}}

\def\bbar{\overline{b}}
\def\cbar{\overline{c}}

\def\ebar{\overline{e}}

\def\hbar{\overline{h}}

\def\wbar{\overline{w}}
\def\xbar{\overline{x}}
\def\ybar{\overline{y}}
\def\zbar{\overline{z}}
\def\Mbar{\overline{M}}

\def\acl{{\rm acl}}

\def\tp{{\rm tp}}
\def\stp{{\rm stp}}

\def\tr/{trivial}
\def\nt/{non--trivial}
\def\st/{strong type}

\def\TV/{Tarski--Vaught}

\def\sc/{sound construction}
\def\ac/{atomic construction}

\def\fal/{functional}
\def\upl/{unique parallel lines}
\def\chp/{categorical in a higher power}

\def\conc{{\char'136}}
\def\abar{\bar{a}}
\def\bbar{\bar{b}}
\def\cbar{\bar{c}}

\def\ebar{\bar{e}}
\def\xbar{\bar{x}}
\def\ybar{\bar{y}}
\def\zbar{\bar{z}}
\def\phi{\varphi}

\def\A{{\cal A}}
\def\E{{\cal E}}
\def\C{{\frak  C}}

\def\MA{{\cal MA}}
\def\tp{{\rm tp}}
\def\stp{{\rm stp}}

\def\acl{{\rm acl}}

\def\bp{{\bf Proof.}\quad}
\def\endproof{\medskip}
\def\<{\langle}
\def\>{\rangle}
\newtheorem{Theorem}{Theorem}[section]
\newtheorem{Proposition}[Theorem]{Proposition}
\newtheorem{Definition}[Theorem]{Definition}

\newtheorem{Example}[Theorem]{Example}
\newtheorem{Lemma}[Theorem]{Lemma}
\newtheorem{Corollary}[Theorem]{Corollary}

\begin{document}
\maketitle
\begin{abstract}
We introduce the notions of a mutually algebraic structures and theories and prove many equivalents.
A theory $T$ is mutually algebraic if and only if it is weakly minimal and  trivial
if and only if no model $M$ of $T$ has an expansion $(M,A)$ by a unary predicate with the finite cover property.
We show that every structure has a maximal mutually algebraic reduct, and give a strong
structure theorem for the class of elementary extensions of a fixed mutually algebraic structure.
\end{abstract}

\section{Introduction}

This paper is written with two objectives in mind.  On one hand, it
is a continuation of \cite{MCL}, where a strong quantifier elimination
theorem was proved for elementary diagrams of models of a weakly minimal,
trivial theory.  Here, we show that the crucial notion of mutual algebraicity
of a formula (see Definition~\ref{ma}) has meaning in arbitrary structures,
and in fact describes a specific reduct of any structure.  
As well, Theorem~\ref{equiv} reverses the argument in \cite{MCL}.
The  quantifier elimination result described there
can only occur as the elementary diagram of a weakly minimal, trivial theory.

On the other hand, there has been a large body of research about
whether an expansion $(M,A)$ of a given stable structure $M$ by a unary
predicate $A$ remains stable.  Sufficient conditions abound, but the general
question remains open.  Here, also with Theorem~\ref{equiv}, we characterize
those structures $M$ with the property that every unary expansion
$(M,A)$ satisfies the non-finite cover property (nfcp), which is a strengthening
of stability.  

The motivation for this came from the author's reading \cite{BB}, where 
Baldwin and
Baizhanov showed that a non-trivial,
strongly minimal structure $M$ has an unstable expansion $(M,A)$.
Thanks are due to  John Baldwin for a careful reading of this paper, and for pointing out
that an alternate treatment of a portion of Section~4 appears in Section~6 of
\cite{BSh}.

\section{The mutually algebraic reduct of a structure}

We begin by recalling the definition of a mutually algebraic formula.
This notion was introduced by Dolich, Raichev, and the author in \cite{DLR}
and further developed in \cite{MCL}.  However, in both of those papers, the ambient
theory was assumed to have the non-finite cover property (nfcp).  Here, we define the notions
without any ambient assumptions.  We begin by formulating the notion of a mutually algebraic set.

\begin{Definition} \label{maset}  {\em  
Given an arbitrary set $A$ and an integer $n\ge 1$, 
a {\em proper partition of $n$\/} is a partition 
$X\sqcup Y=\{1,\dots,n\}$
where $X,Y$ are disjoint and each is non-empty.
Given such a partition, $\pi_Y$ denotes the projection of $A^n$ onto
the coordinates in $Y$.  

A subset $B\subseteq A^n$
is {\em mutually algebraic\/} if there is a number $K$ so that for any
proper partition of the coordinates $X\sqcup Y=\{1,\dots, n\}$,
the projection $\pi_Y$ restricted to $B$ is at most $K$-to-1.
That is, 
$|\pi_Y^{-1}(\bbar_Y)\cap B|\le K$
for any $\bbar_Y\in \pi_Y(B)$.
}
\end{Definition}

As special cases, note that if either $A$ is finite or $B$ is empty, then
$B$ is mutually algebraic.  Furthermore, for any set $A$, every subset 
$B\subseteq A^1$ is mutually algebraic as there are no proper partitions
of a one element set.

\begin{Definition}   \label{ma}
{\em Let $M$ denote any $L$-structure.  
An $L(M)$-formula $\phi(\zbar)$ is {\em mutually algebraic\/}
if $\phi(M):=\{\abar\in M^{\lg(\zbar)}:M\models\phi(\abar)\}$
is a mutually algebraic subset of $M^{\lg(\zbar)}$.
We let $\MA(M)$ denote the set of all mutually algebraic $L(M)$-formulas.  When $M$ is understood,
we simply write $\MA$.

}
\end{Definition}

To clarify this concept and to set notation, 
given a formula $\phi(\zbar)$, a {\em proper partition of
$\zbar$\/} has the form $\zbar=\xbar\conc\ybar$, where 
$\xbar,\ybar$ are disjoint and
$\lg(\xbar),\lg(\ybar)\ge 1$.  
We do not require $\xbar$ be
an initial segment of $\zbar$ 
but to simplify notation, we write it as if it were.
Then, for any $L$-structure $M$, an $L(M)$-formula $\phi(\zbar)$
is mutually algebraic if and only if there is an integer
$K$ so that
$M\models\forall\ybar\exists^{\le K}\xbar\phi(\xbar,\ybar)$ for every proper partition
$\xbar\conc\ybar$ of $\zbar$.

The reader is cautioned that whether a formula $\phi(\zbar)$ is mutually algebraic or not
depends on the choice of free variables.  In particular, mutual algebraicity is
{\bf not} preserved under adjunction of dummy variables.  The special cases
mentioned above imply that if $M$ is finite, then every $L(M)$-formula
is in $\MA(M)$, and for an arbitrary $M$, every inconsistent formula and
every $L(M)$-formula $\phi(z)$ with
exactly one free variable symbol is mutually algebraic.
Our first easy Lemma gives a semantic interpretation to this notion when $\lg(\zbar)\ge 2$:

\begin{Lemma} \label{semantic}
Let $M$ be any $L$-structure.  The following are equivalent
for any $L(M)$-formula $\phi(\zbar)$ with $\lg(\zbar)\ge 2$:
\begin{enumerate}
\item  $\phi(\zbar)\in\MA(M)$;
\item  There is an integer $K$ so that $M\models\forall x\exists^{\le K}\ybar\phi(x,\ybar)$
for all partitions $\zbar=x\conc\ybar$ with $\lg(x)=1$;
\item  For all $N\succeq M$, for all $\ebar\in N^{\lg(\zbar)}$ realizing $\phi$,
and for all $e\in\ebar$, 
$\ebar\subseteq\acl(M\cup\{e\})$ 
(i.e, every $e'\in\ebar$ is in $\acl(M\cup\{e\})$.
\end{enumerate}
\end{Lemma}

\bp  $(1)\Rightarrow(2)$ is immediate.  

$(2)\Rightarrow(3)$  Fix any $N\succeq M$ and assume $N\models\phi(\ebar)$.  Fix any variable
symbol $x\in\zbar$ and let $e$ be the corresponding element of $\ebar$.
By elementarity, $N\models\exists^{\le K}\ybar \phi(e,\ybar)$, 
so $\ebar\subseteq\acl(M\cup\{e\})$.

$(3)\Rightarrow(1)$  If (1) fails, then for some  proper partition $\zbar=\xbar\conc\ybar$ we have
$M\models\exists\ybar\exists^{\ge r}\xbar\phi(\xbar,\ybar)$.
Thus, by compactness, there is $N\succeq M$ and $\bbar$ from $N$ such that
$N\models\exists^{\ge r}\xbar\phi(\xbar,\bbar)$ for each $r\in\omega$.
By compactness again, there is $N^*\succeq N$ and $\abar\in (N^*)^{\lg(\xbar)}$
such that $\abar\not\subseteq\acl(M\cup\bbar)$, contradicting (3).
\endproof

The following Lemma indicates some of the closure properties of the set $\MA$.
In what follows, when we write $\phi(\xbar,\ybar)\in\MA$, we mean that 
$\xbar$ and $\ybar$ are disjoint sets of variable symbols and
$\phi(\zbar)\in\MA$ where $\zbar=\bar\conc\ybar$,
but that we are concentrating on a specific
proper partition  of $\phi(\zbar)$.

\begin{Lemma}  \label{closure}
Let $M$ be any structure in any language $L$.
\begin{enumerate}
\item  If $\phi(\zbar)\in\MA$, then $\phi(\sigma(\zbar))\in\MA$ for any permutation $\sigma$ of the
variable symbols;
\item  If $\phi(\xbar,\ybar)\in\MA$ and $\abar\in M^{\lg(\ybar)}$, 
then both  $\exists\ybar\phi(\xbar,\ybar)$ and $\phi(\xbar,\abar)\in\MA$;
\item  If $\phi(\zbar)\vdash\psi(\zbar)$ and $\psi(\zbar)\in\MA$, then $\phi(\zbar)\in\MA$;
\item  If $\{\phi_i(\zbar_i):i<k\}\subseteq \MA$, and 
there is some variable $x$ common to every $\zbar_i$,
then $\psi(\wbar):=\bigwedge_{i<k} \phi_i(\zbar_i)\in\MA$, where 
$\wbar=\bigcup_{i<k}\zbar_i$;
\item  If $\phi(\xbar,\ybar)\in\MA$ and $r\in\omega$, then $\theta_r(\ybar):=\exists^{\ge r}\xbar\phi(\xbar,\ybar)\in \MA$.
\end{enumerate}
\end{Lemma}

\bp  The verification of (1), (2), and (3) are immediate.
Concerning (4), we apply Lemma~\ref{semantic}.  Fix $N\succeq M$ and
$\ebar$ such that $N\models\psi(\ebar)$.  Let $x$ denote a variable symbol that
appears in every $\zbar_i$
and let
$e_x$ denote the element of $\ebar$ corresponding
to $x$.   
Similarly, for each $i<k$ let $\ebar_i$ be the subsequence corresponding to $\zbar_i$.
As each $\phi_i(\zbar_i)\in\MA$, $e_x\in\acl(M\cup\{e\})$ for every $e\in\ebar_i$, so
$e_x\in\acl(M\cup\{e\})$ for every $e\in\ebar$.
But also, $e\in\acl(M\cup\{e_x\})$ for every $e\in\ebar$.  Thus, by the transitivity of algebraic closure,
$e\in\acl(M\cup\{e'\})$ for all pairs $e,e'\in\ebar$.  So $\psi(\wbar)\in\MA$ by Lemma~\ref{semantic}.

To establish (5), let $\{\xbar_i:i<r\}$ be disjoint sequences of variable symbols, each disjoint from $\ybar$.
Then $\theta_r(\ybar)$ is equivalent to
$$\exists \xbar_0\exists\xbar_1\dots\exists\xbar_{r-1}\left(\bigwedge_{i<r}\phi(\xbar_i,\ybar)\wedge\bigwedge_{i<j<r} \xbar_i\neq\xbar_j\right)$$
That this formula is in $\MA$ follows by successively applying Clauses (4), (3), and (2).
\endproof

\begin{Definition} {\em For any $L$-structure $M$, let $M_M$ denote the canonical expansion of $M$
to an $L(M)$-structure
formed by adding a constant symbol $c_a$ for each $a\in M$. 
We  let $\MA^*(M)$ denote the set of all $L(M)$-formulas
that are $Th(M_M)$-equivalent to a boolean combination
of formulas from $\MA(M)$.
When $M$ is understood, we simply write $\MA^*$.
}
\end{Definition}

Whereas the definition of $\MA$ was rather fussy, membership in $\MA^*$ is more relaxed, mostly owing
to the fact that $\MA^*$ is closed under adjunction of dummy variables.
Indeed, we will see with Proposition~\ref{reduct} below, for any structure $M$,
$\MA^*(M)$ specifies a reduct of  the canonical expansion $M_M$.

\begin{Lemma}  \label{starclosure}
Let $M$ denote any $L$-structure.
\begin{enumerate}
\item  $\MA^*$ is closed under boolean combinations;
\item  $\MA^*$ is closed under adjunction of dummy variables, i.e., if
$\phi(\zbar)\in\MA^*$ then $\phi(x,\zbar)\in\MA^*$;
\item  For each $k\ge 1$, if $\{\phi_i(x,\ybar_i):i<k\}\subseteq\MA$ and $r\in\omega$,
then each of $\exists^{=r}x\bigvee_{i<k}\phi_i(x,\ybar_i)$,
$\exists^{\le r}x\bigvee_{i<k}\phi_i(x,\ybar_i)$, and
$\exists^{\ge r}x\bigvee_{i<k}\phi_i(x,\ybar_i)$ are in $\MA^*$.
\end{enumerate}
\end{Lemma}

\bp  The proof of (1) is immediate.  For (2), note that $\psi(x):=`x=x$' is in $\MA$, hence in $\MA^*$,
but $\phi(x,\zbar)$ is equivalent to $\phi(\zbar)\wedge\psi(x)$.
The verification of (3) is more substantial.  We argue by induction on $k$ that
for every $r\in\omega$,
$\exists^{=r}x\bigvee_{i<k}\phi_i(x,\ybar_i)\in\MA^*$ for every
$k$-element subset $\{\phi_i(x,\ybar_i):i<k\}$ from $\MA$.  This suffices, as $\MA^*$ is closed under boolean
combinations and the trivial facts that $\exists^{\le r}x\theta$ is equivalent to $\bigvee_{s\le r}\exists^{=s}x\theta$ and
$\exists^{\ge r}x\theta$ is equivalent to $\neg\exists^{\le r-1}x\theta$.

To handle the case when $k=1$, fix any  $\phi(x,\ybar)\in\MA$ and any $r\in\omega$.
By Lemma~\ref{closure}(5), both  $\exists^{\ge r}x\phi(x,\ybar)\in\MA$ and 
$\exists^{\ge r+1}x\phi(x,\ybar)\in\MA$ and  $\exists^{=r}x\phi(x,\ybar)$ is a boolean combination of these.

Next, inductively assume that 
for every $r\in\omega$,
$\exists^{=r}x\bigvee_{i<k}\phi_i(x,\ybar_i)\in\MA^*$ for every
$k$-element subset $\{\phi_i(x,\ybar_i):i<k\}$ from $\MA$.  
Choose any $(k+1)$-element subset 
$\{\phi_i(x,\ybar_i):i\le k\}$ from $\MA$ and choose any $r\in\omega$.
As notation, let $\psi(x,\wbar):=\bigvee_{i<k}\phi_i(x,\ybar_i)$.
By the inclusion/exclusion principle of integers,
the formula $\exists^{=r}x\bigvee_{i\le k}\phi(x,\ybar_i)$, which is equivalent to 
$\exists^{=r}x(\psi(x,\wbar)\vee\phi_k(x,\ybar_k))$, is equivalent to
$$\bigvee_{\stackrel {\scriptstyle a,b\le r}{a+b-c=r}} \bigg(\exists^{=a}x\psi(x,\wbar)\ \wedge\ \exists^{=b}x\phi_k(x,\ybar_k)\ \wedge\
\exists^{=c}x[\psi(x,\wbar)\wedge\phi_k(x,\ybar_k)]\bigg)$$
By the inductive hypothesis $\exists^{=a}x\psi(x,\wbar)\in\MA^*$ and $\exists^{=b}x\phi_k(x,\ybar_k)\in\MA^*$
by the case $k=1$.  Also, note that $\psi(x,\wbar)\wedge\phi_k(x,\ybar_k)$ is equivalent to
$\bigvee_{i<k}\delta_i(x,\ybar_i,\ybar_k)$, where each $\delta_i(x,\ybar_i,\ybar_k):=\phi_i(x,\ybar_i)\wedge\phi_k(x,\ybar_k)$
is in $\MA$ by Lemma~\ref{closure}(4).  Thus, by applying the inductive hypothesis to this $k$-element
subset from $\MA$, we conclude that 
$\exists^{=c}x(\psi(x,\wbar)\wedge\phi_k(x,\ybar_k))\in\MA^*$, completing the proof.
\endproof

\begin{Proposition}  \label{reduct}
For any structure $M$, the set $\MA^*(M)$ is closed under existential quantification.
Thus, the structure with universe $M$, together with the definable sets $MA^*(M)$, is a reduct of the 
canonical expansion $M_M$.
\end{Proposition}

\bp The second sentence follows from the first, since $\MA^*$ is a set of $L(M)$-formulas
closed under boolean combinations.  To establish the first sentence, there are two cases.
First, if the structure $M$ is finite, then every $L(M)$-formula $\phi(\zbar)\in\MA$,
so $\MA^*$ is precisely the elementary diagram of $M$ and there is nothing to prove.
So assume that $M$ is infinite.

Choose $\phi(x,\ybar)\in\MA^*$ and we argue that $\exists x\phi(x,\ybar)$ is equivalent to
a formula in $\MA^*$.  By writing $\phi$ in Disjunctive Normal Form and noting that
disjunction commutes with existential quantification, we may assume that $\phi(x,\ybar)$
has the form $$\bigwedge_{i<k} \beta_i(x,\ybar_i)\wedge \bigwedge_{j<m}\neg\gamma_j(x,\ybar_j)$$
where each $\beta_i$ and $\gamma_j$ are in $\MA$ and the variable $x$ occurs in each of these subformulas.
By Lemma~\ref{closure}(4), if $k\ge 1$, then $\bigwedge_{i<k} \beta_i(x,\ybar_i)\in\MA$, so we may assume
there is at most one $\beta$.  
If there is no $\beta$, then since the model $M$ is infinite,
then for any choice of $\ybar$, $\exists x\phi(x,\ybar)$ always holds.
Thus, we assume that there is exactly one $\beta$, i.e., that $\phi(x,\ybar)$
has the form $\beta(x,\ybar^*)\wedge\bigwedge_{j<m}\neg\gamma_j(x,\ybar_j)$, 
where $\ybar^*$ and each $\ybar_j$ are subsequences of $\ybar$, 
and both $\beta$ and each $\gamma_j$
are from $\MA$.  

We first consider the case where  $\ybar^*$ is empty.
In this case, we may additionally assume that no $\ybar_j$ is empty,
since we could replace $\beta(x)$ by $\beta(x)\wedge\neg\gamma_j(x)$.
Thus, for any choice of $\ybar$, the solution set of
$\bigwedge_{j<m}\neg\gamma_j(x,\ybar_j)$ is a cofinite subset of
$M$.
We have two subcases:
On one hand,
if $\beta(x)$ were algebraic, 
then every solution to $\beta$ lies in $M$, 
hence
$\phi(x,\ybar)$ would be equivalent to 
$\bigvee_{m\in\beta(M)}\phi(m,\ybar)$, which would be in $\MA^*$ by
Lemma~\ref{closure}(2).  
On the other hand, if $\beta(x)$ were non-algebraic, 
then $\beta(x)$ would have infinitely many solutions in $M$,
so $\phi(x,\ybar)$ would have a solution in $M$ for any choice of $\ybar$.
Thus, $\exists x\phi(x,\ybar)$ would always hold.

Finally, assume that $\ybar^*\neq\emptyset$.  
By the definition of mutual algebraicity, there is an integer $K$ so
that $M\models\forall\ybar^*\exists^{\le K} x\beta(x,\ybar^*)$.
For each $j<m$, let $\theta_j(x,\ybar^*,\ybar_j):=\beta(x,\ybar^*)\wedge\gamma_j(x,\ybar_j)$.
By Lemma~\ref{closure}(4), each $\theta_j(x,\ybar^*,\ybar_j)\in\MA$.
Thus,  the formula $\exists x\phi(x,\ybar)$ is equivalent to
$$\bigvee_{r\le K} \left(\exists^{=r}x\beta(x,\ybar^*)\wedge\exists^{<r}x\bigvee_{j<m}\theta_j(x,\ybar^*,\ybar_j)\right)$$
which is in $\MA^*$ by Lemma~\ref{starclosure}.
\endproof

The previous Proposition inspires the following two definitions:

\begin{Definition}  {\em  A structure $M$ is {\em mutually algebraic\/}
if every $L(M)$-formula is in $\MA^*(M)$.
}
\end{Definition}

\begin{Definition}  \label{maxreduct}
 {\em  Let $M$ be any structure.  The {\em mutually algebraic reduct of $M_M$}
is the structure with the same universe as $M$, and whose definable sets are precisely
$\MA^*(M)$.
}
\end{Definition}

Proposition~\ref{reduct} immediately implies that the mutually algebraic
reduct of a structure $M$ is a mutually algebraic structure.

\begin{Lemma} \label{elementaryequivalence}
Mutual algebraicity of structures
is preserved under elementary equivalence.
\end{Lemma}

\bp Suppose that $M$ is a mutually algebraic structure and that $N$ is 
elementarily equivalent to $M$.  
It suffices to show that $\phi(\xbar,\hbar)\in\MA^*(N)$
for any $L$-formula $\phi(\xbar,\ybar)$ (with $\xbar$ and $\ybar$ disjoint and
there are no hidden parameters)
and any $\hbar\in N^{\lg(\ybar)}$.
Given this data, let $\zbar=\xbar\conc\ybar$ and consider the $L$-formula $\phi(\zbar)$.
As $M$ is mutually algebraic, $\phi(\zbar)\in\MA^*(M)$, so there are (finitely many) $L$-formulas
$\delta_i(\zbar,\wbar_i)$ and $\ebar_i$ from $M$ so that (1) $\phi(\zbar)$ is
$Th(M_M)$-equivalent to a boolean combination $\theta(\zbar,\ebar^*)$ of the
$\delta_i(\zbar,\ebar_i)$ ($\ebar^*$ denotes the concatenation of the $\ebar_i$'s);
and (2) There is a number $K$ so that each of the formulas $\delta_i(\zbar,\ebar_i)$ 
satisfy $M\models\forall\ybar'\exists^{\le K}\xbar'\delta_i(\xbar',\ybar',\ebar_i)$ for every
proper partition $\zbar=\xbar'\conc\ybar'$.
Thus, by quantifying out the $\ebar^*$, there is an $L$-sentence $\sigma$ asserting 
that
$$\exists\wbar^*\bigg(\forall\zbar [\phi(\zbar)\leftrightarrow\theta(\zbar,\wbar^*)]
\ \wedge\ \hbox{`each $\delta_i(\zbar,\wbar_i)$ is $K$-mutually algebraic'}\bigg)$$
As $M\models \sigma$, so does $N$.  Choose $\cbar^*$ from $N$ so that
$\phi(\zbar)$ is $Th(N_N)$-equivalent to $\theta(\zbar,\cbar^*)$ and $\theta(\zbar,\cbar^*)$
is equivalent to a boolean combination of $\delta_i(\zbar,\cbar_i)\in\MA(N)$, 
where each $\cbar_i$ is 
the corresponding subsequence of $\cbar^*$.   Finally, rewrite $\zbar$ as $(\xbar,\ybar)$
and substitute $\hbar$ for $\ybar$.  By Lemma~\ref{closure}(2), each of the formulas
$\delta_i(\xbar,\hbar,\cbar_i)\in\MA(N)$ and $\phi(\xbar,\hbar)$ is $Th(N_N)$-equivalent
to the boolean combination $\theta(\xbar,\hbar,\cbar^*)$.  Thus, $\phi(\xbar,\hbar)\in\MA^*(N)$,
as required.
\endproof

The following Lemma is folklore, but a proof is included for the convenience
of the reader.  Recall that a partitioned formula $\phi(\xbar,\ybar)$ does not have the finite
cover property (i.e., has nfcp) with respect to a theory $T$
if there is a number $k$ so that for all sets $\{\cbar_i:i\in I\}$,
the type $\Gamma:=\{\phi(\xbar,\cbar_i):i\in I\}$ is consistent with $T$ 
whenever every $k$-element subset
of $\Gamma$ is consistent with $T$.

\begin{Lemma}  \label{cases}
Let  $M$ be  any structure, and let $\phi(\xbar,\ybar)$ be any partitioned
$L(M)$-formula.
If, for some integer $K$, either 
$M\models\forall\ybar\exists^{< K}\xbar\phi(\xbar,\ybar)$, or 
$M\models\forall\ybar\exists^{< K}\xbar\neg\phi(\xbar,\ybar)$,
then $\phi(\xbar,\ybar)$ does not have the finite cover property
with respect to $Th(M_M)$.

\end{Lemma}

\bp If $M$ is finite, then every partitioned formula $\phi(\xbar,\ybar)$ has nfcp for trivial reasons,
so assume that $M$ is infinite.  First, assume that
$M\models\forall\ybar\exists^{< K}\xbar\phi(\xbar,\ybar)$.
Choose tuples $\{\cbar_i:i\in I\}$ from some elementary extension of $M$
and assume that the type $\Gamma:=\{\phi(\xbar,\cbar_i):i\in I\}$ is inconsistent.
It suffices to find a subtype of at most $K$ elements that is inconsistent as well.
Choose a maximal sequence $\<i_j:j\le n\>$ from $I$ such that $i_0\in I$ is arbitrary
and for each $1\le m\le n$,
$$\models \exists \xbar \left(\bigwedge_{j<m}
\phi(\xbar,\cbar_{i_j})\wedge\neg\phi(\xbar,\cbar_{i_m})\right)$$
By our hypotheses on $\phi(\xbar,\cbar_{i_0})$, $n\le K$.  But now, if
$\bigwedge_{j\le n}\phi(\xbar,\cbar_{i_j})$ were consistent but 
$\Gamma$ were not, we would
contradict the maximality of the sequence.

In the other case, as $M$ is infinite, every partial type of the form 
$\{\phi(\xbar,\cbar_i):i\in I\}$ is consistent, so the nfcp of $\phi(\xbar,\ybar)$ is vacuously true.
\endproof

\begin{Proposition}  \label{nfcp}
For any structure $M$, the theory of the mutually algebraic reduct of $M$ has nfcp.
\end{Proposition}

\bp  
By the equivalence of (1) and $\forall m(2)_m$ in Theorem~II~4.4 of \cite{Shc} (whose proof does not use
stability) it suffices to show that no partitioned formula of the form $\phi(x,\ybar)\in\MA^*$
with $\lg(x)=1$ has the finite cover property.

Consider any formula $\theta(x,\ybar)$ of the form
$$\bigwedge_{i<k} \beta_i(\zbar_i)\wedge\bigwedge_{j<m} \neg\gamma_j(\zbar_j)$$
with each $\beta_i$ and $\gamma_j$ from $\MA$.
First, if the variable $x$ occurs in any $\beta_i$, then it follows that
there is a number $K$ so that $M\models\forall\ybar\exists^{< K}x\theta(x,\ybar)$.
Second, if $x$ does not occur in any $\beta_i$, then there is a number $K$ so
that there is a number $K$ so that $M\models\forall\ybar\exists^{< K}x\neg\theta(x,\ybar)$.
But, any formula $\phi(x,\ybar)\in\MA^*$ is a finite disjunction of formulas
$\theta(x,\ybar)$ described above.  It follows that for some $K$,
either
$M\models\forall\ybar\exists^{< K}x\phi(x,\ybar)$
or
there is a number $K$ so that $M\models\forall\ybar\exists^{< K}x\phi(x,\ybar)$ holds.
Thus, $\phi(x,\ybar)$ has the nfcp by Lemma~\ref{cases}.

\section{Characterizing theories of mutually algebraic structures}

We begin with two definitions indicating that the forking behavior of 1-types (types with a
single free variable) is particularly simple.

\begin{Definition} {\em  A complete, stable theory with an infinite model
is {\em weakly minimal\/} if
every forking extension of a 1-type is algebraic (equivalently if $R^\infty(x=x)=1$)
and is {\em trivial} if there do not
exist a set $D$ and three elements $\{a,b,c\}$ that are dependent, but pairwise independent over $D$. 
A type $p\in S(D)$ is trivial if there do not exist a set $\{a,b,c\}$ of realizations of
$p$ that are dependent, but pairwise independent over $D$. }
\end{Definition}

It is well known that a weakly minimal theory is trivial if and only if
every minimal type is trivial.
The following Lemma  generalizes the analogous result for non-trivial,
strongly minimal theories that was proved by Baldwin and Baizhanov in \cite{BB}.

\begin{Lemma} \label{wmtrivial}
 If $T$ is weakly minimal and non-trivial, then
there is a model $M$ of $T$ and a subset $A\subseteq M$ such
that $(M,A)$ is unstable.
\end{Lemma}

\bp  Among all minimal types $p\in S(D)$ and formulas
$\phi(z,xy)$ over $D$ that contain a dependent, but pairwise
independent triple $\{a,b,c\}$ of realizations of $p$, with the dependency 
witnessed by the algebraic formula $\phi(z,ab)\in\tp(c/Dab)$, choose one
with the multiplicity of $\phi(z,ab)$ as small as possible.  It follows
from this multiplicity condition that $\acl(D\cup\{a\})\cup\acl(D\cup\{b\})$
does not contain any realizations of $\phi(z,ab)$.

Fix $p\in S(D)$ and $\phi(z,xy)$ as above, and let
$M$ be a sufficiently saturated model containing $D$. To ease notation, we may
assume $D=\emptyset$.  Let $\<(a_i,b_i):i\in\omega\>$ be a Morley
sequence in $p^{(2)}$.  That is, $\{a_i:i\in\omega\}\cup\{b_j:j\in\omega\}$
is an independent set of realizations of $p$.  For each pair $(i,j)\in\omega^2$,
choose $c_{i,j}\in p(M)$ realizing $\phi(z,a_ib_j)$.
Let $A=\{c_{i,j}:i\le j<\omega\}$. We argue that the $L_P$-formula
$\Phi(x,y):=\exists z (P(z)\wedge\phi(z,xy))$ has the order property in
$(M,A)$.  

To see this, it is clear that the element $c_{i,j}$ witnesses $\Phi(a_i,b_j)$
whenever $i\le j$.  On the other hand, suppose 
some $c_{k,\ell}$ witnessed $\phi(x,a_i,b_j)$.  We argue that we must
have $k=i$ and $\ell=j$:  If neither equality held, then we would
have $c_{k,\ell}$ forking with both sets $\{a_i,b_j\}$ and $\{a_k,b_{\ell}\}$.
This is impossible, as the doubletons are independent from each other
and the type $p$ is  minimal, hence regular, hence of weight one.
Similarly, suppose that $k=i$ but $\ell\neq j$.  Then, working over $a_i$,
$c_{i,\ell}$ is not algebraic over $a_i$, so $\tp(c_{i,\ell}/a_i)$ is parallel
to $p$, hence is also regular, so of weight one.  But, working over $a_i$,
$c_{i,\ell}$ forks with each of $b_j$ and $b_{\ell}$, which are independent
over $a_i$.  The case where $j=\ell$ is symmetric, completing the proof.
\endproof

In what follows, a {\em mutually algebraic expansion\/} of a structure $M$
is an expansion formed by adding arbitrarily many new relation symbols $R_i$,
whose interpretation is a mutually algebraic subset of $M^{{\rm arity}(R_i)}$
(see Definition~\ref{maset}).
In the Theorem that follows, we do not require that the theory $T$ be complete.

\begin{Theorem}  \label{equiv}
The following are equivalent for any theory $T$:
\begin{enumerate}
\item  Every model of $T$ is a mutually algebraic structure;
\item  Every mutually algebraic expansion of every model of $T$
is a mutually algebraic structure;
\item  $Th((M,A))$ has the nfcp for every  $M\models T$ and every expansion $(M,A)$ by a unary predicate;
\item  Every complete extension of $T$ having an infinite model is weakly minimal and trivial.
\end{enumerate}
\end{Theorem}

\bp  $(1)\Rightarrow(2)$ Fix $M\models T$ and let $\Mbar=(M,R_i)_{i\in I}$ be any expansion of
$M$, where each $R_i$ is a $k(i)$-ary relation symbol whose interpretation in $\Mbar$ is 
a mutually algebraic subset $B_i\subseteq M^{k(i)}$.  By definition, the $\Mbar$-definable subsets
are the smallest class of subsets of $M^{\ell}$ for various $\ell$ that contain every $M$-definable
set and every $B_i$ and are closed under boolean combinations and projections.
As $M$ is mutually algebraic, every $M$-definable set is a boolean combination of mutually algebraic sets.  So
$\MA^*(\Mbar)$ contains every $M$-definable set and each of the sets $B_i$.  Additionally, $\MA^*(\Mbar)$
is closed under boolean combinations and projections.  Thus, every $\Mbar$-definable set
is in $\MA^*(\Mbar)$, so $\Mbar$ is a mutually algebraic structure.

$(2)\Rightarrow(3)$  Fix any $M\models T$
and any expansion $\Mbar=(M,A)$ by a unary predicate.
As every subset of $M^1$ is mutually algebraic, it follows from (2) that $\Mbar$ is a mutually algebraic
structure, i.e., every $\Mbar$-definable set is in $\MA^*(\Mbar)$.  Thus, every partitioned
$\Mbar$-definable formula $\phi(\xbar,\ybar)$ has nfcp by Proposition~\ref{nfcp}.
That is, the elementary diagram of $\Mbar$ and hence the theory of $\Mbar$ has nfcp.

$(3)\Rightarrow(4)$  Suppose $T$ satisfies (3).
Fix any complete extension $T'$ of $T$ with an infinite model.  
As the nfcp implies stability, $T'$ must be
stable.  Fix a sufficiently saturated model $M$ of $T'$.
As $T'$ is stable, if it were not weakly minimal then we could choose an
element $a$ and a tuple $\bbar$ from $M$ such that $\tp(a/\bbar)$ forks over the empty set,
but $a$ is not algebraic over $\bbar$.  Let $\phi(x,\ybar)$ be chosen so
that $\phi(x,\bbar)\in\tp(a/\bbar)$ witnesses the forking.  As $M$ is sufficiently saturated, we
can find a Morley sequence $\<\bbar_i:i\in\omega\>$ in $\stp(\bbar)$ inside $M$.
As $T'$ is stable, $\{\bbar_i:i\in\omega\}$ is an indiscernible set and there is a number $k$
so that every element $a^*\in M$ is contained in at most $k$ of the sets $D_i:=\phi(M,\bbar_i)$.
As each $D_i$ is infinite, we can construct a subset $A$ of $M$
such that each $c\in A$ is contained in exactly one of the sets $D_i$, and for each
$i$, $|A\cap D_i|=i$.  Then the theory of the expansion $(M,A)$, where the new unary predicate
symbol $P$ is interpreted as $A$, has the finite cover property
as witnessed by the $L_P$-formula $\Psi(x,\ybar z):=P(x)\wedge\phi(x,\ybar)\wedge x\neq z$. 
Thus, $T$ must be weakly minimal.  That $T'$ must be trivial as well follows from 
Lemma~\ref{wmtrivial} and the fact that instability implies an instance of the finite
cover property.

$(4)\Rightarrow(1)$  This is the content of Theorem~4.2 of \cite{MCL}.  In fact, there it is shown
that every $M$-definable formula is a boolean combination of mutually algebraic formulas of a
very special form.
\endproof

\begin{Corollary}  Let $M$ be any infinite structure.  The mutually algebraic
reduct of $M$ described in Definition~\ref{maxreduct}
is the maximal weakly minimal, trivial reduct of $M$.
\end{Corollary}

\bp  The mutually algebraic reduct of $M$ is a mutually algebraic structure,
so it has a weakly minimal, trivial theory.  Conversely, if any reduct of $M$
has a weakly minimal, trivial theory, then it is a
mutually algebraic structure, hence all of its definable
sets are contained in $\MA^*(M)$.
\endproof

\section{Mutually algebraic structures}

Suppose that $M$ is a mutually algebraic structure in a language $L$. 
We study models of
the elementary diagram of $M$, or equivalently the class of elementary extensions of $M$.
Note that if $M$ is finite, then there are no proper elementary extensions of $M$, which will
render all of the results that follow vacuous.  Because of this, {\bf throughout this section
we additionally assume that $M$ is infinite.}
Thus, we may assume that $M$ is elementarily embedded in a
much larger, saturated `monster model' $\C$.  

By Theorem~\ref{equiv}, $Th(M)$ is weakly minimal and trivial, so the quantifier elimination offered
in \cite{MCL} applies.  Specifically,  let
$$\A(M):=\{\hbox{all quantifier-free mutually algebraic $L(M)$-formulas}\ \alpha(\zbar)\} \ \hbox{and} $$
$$\E(M)=\{\hbox{all $L(M)$-formulas of the form $\exists\xbar\alpha(\xbar,\ybar)$, where}\ \alpha(\xbar,\ybar)\in\A(M)\}$$
and let $\A^*(M)$ (respectively $\E^*(M)$) denote the closure of $\A(M)$ (respectively $\E(M)$) under boolean combinations.
Proposition~4.1 of \cite{MCL} states that every quantifier-free $L(M)$-formula is equivalent to a formula
in $\A^*(M)$, while Theorem~4.2 states that every $L(M)$-formula is equivalent to a formula in $\E^*(M)$.

As $Th(M)$ is weakly minimal, the relation `$a\in\acl_M(B)$' satisfies the axioms of a pre-geometry,
where $\acl_M(B)$ abbreviates $\acl(M\cup B)$.   (Algebraic closures are
always computed with respect to satisfaction in $\C$.)
Thus, the binary relation $a\approx b$ on $\C\setminus M$ defined by
$a\in\acl_M(\{b\})$ is an equivalence relation.   The following easy Lemma
is folklore.

\begin{Lemma}  \label{acl}
Suppose $Th(M)$ is weakly minimal, $M\preceq\C$, and $A$ is any 
algebraically closed set
satisfying $M\subseteq A\subseteq\C$.  Then $A$ is the universe of an
elementary submodel of $\C$.
\end{Lemma}

\bp  The interpretation of any constant symbol is contained in $M$,
and the fact that $A$ is algebraically closed implies that it is closed
under every function symbol in the language.  Thus, $A$ is the universe of
a substructure of $\C$.  To see that this substructure is elementary, by
the Tarski-Vaught criterion it suffices to show that for any $L$-formula
$\phi(x,\ybar)$ and for any $\abar$ from $A$, if $\C\models\exists x\phi(x,\abar)$,
then there is $b\in A$ such that $\C\models\phi(b,\abar)$.
So fix any $\phi(x,\abar)$ and $b\in\C$ such that $\C\models\phi(b,\abar)$.
If $b\in A$, then we are done, so assume $b\not\in A$.  As $A$ is algebraically
closed, this means that $\tp(b/M\abar)$ is not algebraic.  As $Th(M)$ is
weakly minimal and $b$ is a singleton, this implies that $\tp(b/M\abar)$
does not fork over $M$.
But then, by symmetry and finite satisfiability of non-forking over models,
there is $b^*\in M\subseteq A$ such that $\C\models\phi(b^*,\abar)$.
\endproof

Recall that when combined with weak minimality,
triviality implies that for any set $B\subseteq\C$, $\acl_M(B)=\bigcup_{b\in B}\acl_M(\{b\})$.

\begin{Proposition}  \label{model}  Let $M$ be any mutually algebraic $L$-structure.
\begin{enumerate}
\item  If $M\subseteq A\subseteq \C$ and $A$ is an arbitrary union of $\approx$-classes, then
$A$ is an $L$-structure and $M\preceq A\preceq\C$; and
\item  Conversely, if $M\preceq N\preceq\C$ and $B\subseteq N\setminus M$ is a
set of $\approx$-representatives, then $N$ is the disjoint union of the sets $M$ and
$\{\acl_M(\{b\})\setminus M:b\in B\}$.
\end{enumerate}
\end{Proposition}

\bp  (1)  $Th(M)$ is weakly minimal and trivial by Theorem~\ref{equiv}.
By triviality, $A$ must be algebraically closed, so
$A\preceq \C$ by Lemma~\ref{acl}.  That $M\preceq A$
follows immediately from this.

(2)  That the sets are disjoint follows by triviality.  If there were an element $d\in N$ that was
not in any of these sets, then $d$ would be $\approx$-inequivalent to every element of $B$, contradicting
the maximality of $B$.
\endproof

In light of the previous Proposition, it is natural to refer to the sets $\acl_M(\{b\})\setminus M$
as the {\em components} of a given $N\succeq M$.  Each component has size bounded by the number of
$L(M)$-formulas, and one can speak of the type of a fixed enumeration of a component over $M$.
The notion of a component map  records this amount of data.

\begin{Definition} {\em  Suppose that $M$ is a mutually algebraic structure
 and $N_1, N_2$ are both
elementary extensions of $M$.  A {\em component map} $f:N_1\rightarrow N_2$ is a bijection
such that $f|_M=id$ and for each $b\in N_1\setminus M$,
\begin{itemize}
\item  $f$ restricted to $\acl_M(\{b\})$ is elementary and
\item  $f(\acl_M(\{b\}))=\acl_M(\{f(b)\})$ setwise.
\end{itemize}
}
\end{Definition}

\begin{Proposition}  \label{iso}  Suppose that $M$ is mutually algebraic and $N_1,N_2\succeq M$.
Then every component map $f:N_1\rightarrow N_2$ is an isomorphism.
Conversely, every isomorphism $f:N_1\rightarrow N_2$ that is the identity on $M$
is a component map.
\end{Proposition}

\bp 
As every quantifier-free $L(M)$-formula is equivalent to a formula in $\A^*(M)$, it
suffices to show that $f$ preserves every formula $\alpha(\zbar)\in\A(M)$.  Choose any $\abar$ from $N_1$.
Without loss, by Lemma~\ref{closure}(2) and the fact that $f$ fixes $M$ pointwise, we may assume
$\abar$ is disjoint from $M$.  There are now two cases.  First, if $\abar\subseteq\acl_M(\{b\})$ for some
element $b$, then $f(\abar)\subseteq\acl_M(\{b\})$ and $N_1\models\alpha(\abar)$ if and only if
$N_2\models \alpha(f(\abar))$ by the elementarity of $f$ restricted to $\acl_M(\{b\})$.
Second, if $\abar$ intersects at least two components, then $N_1\models\neg\alpha(\abar)$ automatically.
Furthermore, since $f$ maps components onto components, $f(\abar)$ would intersect at
least two components of $N_2$, so $N_2\models\neg\alpha(\abar)$.  Thus, $\alpha$
is preserved in both cases, so $f$ is an isomorphism.
The converse is clear since elementary maps preserve algebraic closure.
\endproof

We close with two examples of how the analog of Proposition~\ref{iso} can fail 
if we work over $\acl(\emptyset)$ instead of a model.

\begin{Example}  Let $L=\{R,S,E\}$, and let $T$ be the theory asserting that $E$ is an equivalence
relation with exactly two classes, both infinite, and 
$R$ is a binary `mating relation' i.e., $R$ is symmetric,
irreflexive, and $\forall x\exists^{=1}y R(x,y)$.  We further require that $R(x,y)\rightarrow \neg E(x,y)$.
Take $S$ to be a 4-ary relation such that $S(x,y,z,w)$ holds if and only if the four elements are distinct,
and each of the relations $R(x,y)$, $R(z,w)$, and $E(x,z)$ hold.  
Then $T$ is complete, mutually algebraic, and
$\acl(\emptyset)=\emptyset$. For any model $N$ of $T$, the decomposition of
$N$ into $R$-mated pairs is a  decomposition  of $N$ into two-element `$\emptyset$-components'
i.e., sets  $A$ satisfying $\acl_\emptyset(A)=A$.
However, in contrast to Proposition~\ref{iso}, there are `$\emptyset$-component maps' $f:N\rightarrow N$, i.e.,
bijections $f:N\rightarrow N$ whose restriction to each two-element $\emptyset$-component is elementary, that
are not automorphisms.
\end{Example}

The second example is from \cite{BLS}.  There, Baldwin, Shelah, and the author exhibit two models
$M, N$ of the theory of infinitely many, binary splitting equivalence relations 
that are not isomorphic in the set-theoretic universe
$V$, but there is a c.c.c.\ extension $V[G]$ of $V$ and $M\cong N$ in $V[G]$.
This theory is also weakly minimal and trivial with $\acl(\emptyset)=\emptyset$.
In fact, this theory has a prime model and every `component' is a singleton.
The complexity exploited by this example involves which strong types over
the empty set are realized in the models $M$ and $N$.

\end{document}